\documentclass[english,11pt,a4paper]{article}
\usepackage[latin9]{inputenc}
\usepackage{amsmath}
\usepackage{amsthm}
\usepackage{makeidx}
\usepackage{latexsym}
\usepackage{stackrel}
\usepackage{babel}
\usepackage{hyperref}
\usepackage{float}
\newtheorem{theorem}{Theorem}
\newtheorem{lemma}{Lemma}
\newtheorem{proposition}{Proposition}

\newtheorem{remark}{Remark}
\newtheorem{definition}{Definition}
\theoremstyle{remark}

\title{Generalized Vandermonde Determinants and Characterization of Divisibility Sequences}
\author{Stefano Barbero\\ \\
	Department of Mathematics, University of Turin\\
	Via Carlo Alberto 10, 10122, Turin, ITALY\\
 stefano.barbero@unito.it}
\date{}

\begin{document}
	\maketitle
	
	\begin{abstract}
	\noindent We present a different proof of the characterization of non--degenerate recurrence sequences, which are also divisibility sequences, given by  Van der Poorten, Bezevin, and Peth{\"o} in their  paper \cite{BPP}. Our proof is based on an interesting  determinant identity related to impulse sequences, arising from the evaluation of a generalized Vandermonde determinant. As a consequence of this new proof we can find a more precise form for the resultant sequence presented in \cite{BPP}, in the general case of non--degenerate divisibility sequences having minimal polynomial with multiple roots.
	\end{abstract}
	\makeatletter{\renewcommand* \@makefnmark{}
\makeatother}
\section{Introduction}
Finding properties for non--degenerate recurrence sequences and also divisibility sequences and determining some kind of deeper structure characterizing  them is a very fascinating research field. 
  The most important attempt to establish their behaviour in an elegant way was presented in the paper of Van der Poorten, Bezevin, and Peth{\"o}  \cite{BPP}, where  they confirm what Ward conjectured in his paper \cite{Ward} about the possibility that every linear divisibility sequence should be a divisor of a resultant sequence. In a field $F$ of characteristic zero, they considered a non--degenerate recurrence sequence  $(a_n)_{n=0}^{+\infty}$, with characteristic polynomial having distinct roots. Using the Hadamard quotient theorem and the theory of exponential polynomials they stated that  if such a sequence  is  a divisibility sequence,  then  there is a resultant sequence $(\bar{a}_n)_{n=0}^{+\infty}$ such that
$$\forall n\geq 0  \quad a_n|\bar{a}_n,$$ 
  where $\bar{a}_n$  has the  shape
\begin{equation}\label{res1}
\bar{a}_n=n^{k}\underset{i}{\prod}\left(\frac{\alpha_{i}^{n}-\beta_{i}^{n}}{\alpha_{i}-\beta_{i}}\right).
\end{equation}

The aim of this paper is to present a proof of this result based on generalized Vandermonde determinants. We start proving an interesting identity concerning non--degenerate impulse sequences and generalized Vandermonde determinants. Then we use it to restate the main result presented in \cite{BPP}, giving a refinement and a more precise form for the $n$--th term of the resultant sequence involved. We deal with the general case of non--degenerate recurrence sequences which are also divisibility sequences and whose minimal polynomial has multiple roots. From now on we work over a field $F$ of characteristic zero. We also remember, once and for all, that we consider a recurrence sequence as non--degenerate if the ratio of two distinct roots of its minimal polynomial is not a root of unity, and obviously all the roots are different from zero.

\section{Impulse sequences and generalized Vandermonde determinants}

We recall the definition of the particular recurrence sequences named \emph{impulse sequences}.

\begin{definition}\label{impulse}
We define the \emph{impulse sequences} of order $r$ as the non--degenerate linear recurrence sequences $\left(X_{n}^{\left(k\right)}\right)_{n=0}^{+\infty}$,
$k=0,\ldots,r-1$,
starting with the initial conditions $X_{j}^{\left(k\right)}=\delta_{jk}$, $j=0,\ldots,r-1,$
($\delta_{jk}$ is the usual Kronecker delta),  whose minimal polynomial has  $s$ distinct roots $\alpha_i$, of respective multiplicities $m_l$, with $l=1,\ldots,s$, $\stackrel[l=1]{s}{\sum}m_{l}=r$. 
\end{definition}

In the next theorem we prove a determinant identity involving these sequences, which will allow us to give in the next section an elementary proof of the characterization of divisibility sequences presented in \cite{BPP}. This identity plainly connects impulse sequences to generalized Vandermonde determinants. During the proof of this Theorem we will use the following lemma based on the results of Flowe and Harris \cite{FH} and exposed as Theorem 21 in the wonderful compendium on determinant calculus written by Krattenthaler \cite{K}.
\begin{lemma}\label{lemfh}
	Let $r$ be a nonnnegative integer, and let $B_m(x)$ denote the $r\times m$ matrix 
	\[
	B_{m}(X)=\left[\begin{array}{cccc}
	1 & 0 & \cdots & 0\\
	X & X & \cdots & X\\
	X^2 & 2X^2 & \cdots & 2^{m-1}X^2\\
	x^3&  3X^3 & \cdots & 3^{m-1}X^3\\
	\cdots & \cdots & \cdots & \cdots\\
	X^{(r-1)} & (r-1)X^{(r-1)} & \cdots & (r-1)^{m-1}X^{(r-1)}
	\end{array}\right]
	\]

	i.e, any next column is formed by applying the operator $X(\frac{d}{dX})$. Given a composition of $r=m_1+m_2+\cdots+m_s$ , there holds
	$$\underset{1\leq i,j\leq r}{\det}
	\left(B_{m_1}(X_1)B_{m_2}(X_2)\cdots B_{m_s}(X_s)	\right)=\stackrel[l=1]{s}{\prod}\stackrel[j=1]{m_{l}-1}{\prod}j!X_{l}^{\left(\begin{array}{c}
		m_{l}\\
		2
		\end{array}\right)}\underset{1\leq i<j\leq s}{\prod}\left(X_{j}-X_{i} \right)^{m_{i}m_{j}}$$
\end{lemma}
\begin{proof}
	For an exhaustive proof we refer the reader to the paper of Flowe and Harris \cite{FH}
\end{proof}

\begin{theorem}\label{detimp}
	Let us consider the $r$ impulse sequences $\left(X_{n}^{\left(k\right)}\right)_{n=0}^{+\infty}$ introduced in Definition \ref{impulse} and the determinant
	
	\begin{equation}\label{vander}
D=	\det\left|\begin{array}{cccc}
	X_{n}^{(1)} & X_{n}^{(2)} & \cdots & X_{n}^{(r-1)}\\
	X_{2n}^{(1)} & X_{2n}^{(2)} & \cdots & X_{2n}^{(r-1)}\\
	\cdots & \cdots & \ddots & \cdots\\
	X_{(r-1)n}^{(1)} & X_{(r-1)n}^{(2)} & \cdots & X_{(r-1)n}^{(r-1)}
	\end{array}\right|.\end{equation}
	Then we have
	\begin{equation}\label{vanderval}
D=	n^{\stackrel[l=1]{s}{\sum}\left(\begin{array}{c}
		m_{l}\\
		2
		\end{array}\right)}\stackrel[l=1]{s}{\prod}\alpha_{l}^{\left(\begin{array}{c}
		m_{l}\\
		2
		\end{array}\right)(n-1)}\underset{1\leq i<j\leq s}{\prod}\left(\frac{\alpha_{j}^{n}-\alpha_{i}^{n}}{\alpha_{j}-\alpha_{i}}\right)^{m_{i}m_{j}}.
		\end{equation}
	
\end{theorem}
\begin{proof}

In order to explicitly evaluate $D$ we point out (see, e.g., the fundamental book on recurrence sequences \cite{EPSW}) that for $k=1,\ldots,r-1$

\begin{equation}
X_{hn}^{(k)}=\stackrel[j=1]{s}{\sum}\stackrel[i=0]{m_{j}-1}{\sum}c_{i,j}^{(k)}(nh)^{i}\alpha_{j}^{nh}\quad h=1,\ldots,r-1\label{eq:xn}
\end{equation}

where, from the initial conditions on $\left(X_{n}^{\left(k\right)}\right)_{n=0}^{+\infty}$,
the coefficients $c_{i,j}^{(k)}$ must satisfy the relations 

\begin{equation}
\stackrel[j=1]{s}{\sum}\stackrel[i=0]{m_{j}-1}{\sum}c_{i,j}^{(k)}t^{i}\alpha_{j}^{t}=\delta_{k,t}\quad t=0,\ldots,r-1\label{eq:ci}
\end{equation}

(with the convention $0^{0}=1$). From these relations (\ref{eq:xn})
and (\ref{eq:ci}) we observe that $D$ is related to the product
between the determinants of the following matrices

\begin{equation}
W\left(n\right)=\left[A_{1}A_{2}\cdots A_{s-1}A_{s}\right]\quad C=\left[\begin{array}{c}
C_{1}\\
C_{2}\\
\vdots\\
C_{s-1}\\
C_{s}
\end{array}\right],\label{eq:matblock}
\end{equation}
 where for $l=1,\ldots,s$ every block $A_{l}$ is an $r\times m_{l}$
matrix and every block $C_{l}$ is an $m_{l}\times r$ matrix, such
that 

\[
A_{l}=\left[\begin{array}{cccc}
1 & 0 & \cdots & 0\\
\alpha_{l}^{n} & n\alpha_{l}^{n} & \cdots & n^{m_{l}-1}\alpha_{l}^{n}\\
\alpha_{l}^{2n} & 2n\alpha_{l}^{2n} & \cdots & (2n)^{m_{l}-1}\alpha_{l}^{2n}\\
\cdots & \cdots & \cdots & \cdots\\
\alpha_{l}^{(r-1)n} & (r-1)n\alpha_{l}^{(r-1)n} & \cdots & [(r-1)n]^{m_{l}-1}\alpha_{l}^{(r-1)n}
\end{array}\right]
\]

\[
C_{l}=\left[\begin{array}{ccccc}
c_{0,l}^{(1)} & c_{0,l}^{(2)} & \cdots & c_{0,l}^{(r-1)} & \delta_{s,l}\\
c_{1,l}^{(1)} & c_{1,l}^{(2)} & \cdots & c_{1,l}^{(r-1)} & 0\\
c_{2,l}^{(1)} & c_{2,l}^{(2)} & \cdots & c_{2,l}^{(r-1)} & 0\\
\cdots & \cdots & \cdots & \cdots & 0\\
c_{m_{l}-1,l}^{(1)} & c_{m_{l}-1,l}^{(2)} & \cdots & c_{m_{l}-1,l}^{(r-1)} & 0
\end{array}\right]
\]
In fact we easily obtain 
\begin{equation}
(-1)^{r+1}D=\det(W\left(n\right))\det(C).\label{eq:reld}
\end{equation}

Moreover from (\ref{eq:ci}), we have

\begin{equation}
W^{*}\left(1\right)C=\left[\begin{array}{cccccc}
1 & 0 & 0 & \cdots & 0 & \alpha_{s}\\
0 & 1 & 0 & \cdots & 0 & \alpha_{s}^{2}\\
0 & 0 & 1 & \cdots & 0 & \alpha_{s}^{3}\\
\cdots & \cdots & \cdots & \ddots & \cdots & \cdots\\
0 & 0 & 0 & \cdots & 1 & \alpha_{s}^{r-1}\\
0 & 0 & 0 & \cdots & 0 & 1
\end{array}\right]\label{eq:wic}
\end{equation}

where

\[
W^{*}\left(1\right)=[A_{1}^{*}A_{2}^{*}\cdots A_{s-1}^{*}A_{s}^{*}]
\]

and every block $A_{l}^{*}$ for $l=1,\ldots,s$ is an $r\times m_l$ matrix of the form

\[
A_{l}^{*}=\left[\begin{array}{cccc}
\alpha_{l} & \alpha_{l} & \cdots & \alpha_{l}\\
\alpha_{l}^{2} & 2\alpha_{l}^{2} & \cdots & (2)^{m_{l}-1}\alpha_{l}^{2}\\
\cdots & \cdots & \cdots & \cdots\\
\alpha_{l}^{(r-1)} & (r-1)\alpha_{l}^{(r-1)} & \cdots & [(r-1)]^{m_{l}-1}\alpha_{l}^{(r-1)}\\
1 & 0 & \cdots & 0
\end{array}\right].
\]
 Clearly from (\ref{eq:wic}) we get

\begin{equation}
\det\left(W^{*}\left(1\right)\right)\det\left(C\right)=(-1)^{r-1}\det\left(W\left(1\right)\right)\det(C)=1.\label{eq:proddetc}
\end{equation}

The last step is to evaluate $\det\left(W\left(n\right)\right)$.
From all the consecutive $m_{l}-1$ columns of $W\left(n\right)$
with first entry equal to 0, we can pick up the factors $n,n^{2},\ldots,n^{m_{l}-1}$
with $l=1,\ldots,s$ . Thanks to Lemma \ref{lemfh}, the determinant of the so--obtained matrix $\stackrel{\sim}{W}\left(n\right)$ satisfies the equality

\[
\det\left(\stackrel{\sim}{W}\left(n\right)\right)=\stackrel[l=1]{s}{\prod}\stackrel[j=1]{m_{l}-1}{\prod}j!\alpha_{l}^{\left(\begin{array}{c}
m_{l}\\
2
\end{array}\right)n}\underset{1\leq i<j\leq s}{\prod}\left(\alpha_{j}^{n}-\alpha_{i}^{n}\right)^{m_{i}m_{j}}
\]
since we only made the substitutions $X_i=\alpha_{i}^n$ for $i=1,\ldots,s$.

Therefore, taking in account the product of the terms picked up in
order to find $\stackrel{\sim}{W}\left(n\right)$, we easily obtain

\[
\det\left(W\left(n\right)\right)=n^{\stackrel[l=1]{s}{\sum}\left(\begin{array}{c}
m_{l}\\
2
\end{array}\right)}\stackrel[l=1]{s}{\prod}\stackrel[j=1]{m_{l}-1}{\prod}j!\alpha_{l}^{\left(\begin{array}{c}
m_{l}\\
2
\end{array}\right)n}\underset{1\leq i<j\leq s}{\prod}\left(\alpha_{j}^{n}-\alpha_{i}^{n}\right)^{m_{i}m_{j}}.
\]
Since we consider non--degenerate recurrence sequences we have $$\forall n\geq 1 \quad \det\left(W(n)\right)\neq0.$$
Now combining (\ref{eq:proddetc}) and (\ref{eq:reld}) we plainly have

\begin{equation}
D=\frac{\det\left(W\left(n\right)\right)}{\det\left(W\left(1\right)\right)}=n^{\stackrel[l=1]{s}{\sum}\left(\begin{array}{c}
	m_{l}\\
	2
	\end{array}\right)}\stackrel[l=1]{s}{\prod}\alpha_{l}^{\left(\begin{array}{c}
	m_{l}\\
	2
	\end{array}\right)(n-1)}\underset{1\leq i<j\leq s}{\prod}\left(\frac{\alpha_{j}^{n}-\alpha_{i}^{n}}{\alpha_{j}-\alpha_{i}}\right)^{m_{i}m_{j}}.
\end{equation}
\end{proof}
\section{Characterizing property of divisibility sequences}
First of all, we give the definition of a non--degenerate linear recurrence sequence which is also a divisibility sequence.
\begin{definition}\label{divseq}
Let us consider a non--degenerate linear recurrence sequence $S=\left(S_{n}\right)_{n=0}^{+\infty}$
of order $r$ with minimal polynomial having $s$ distinct roots $\alpha_{i}$, with respective multiplicities $m_{l}$, where $l=1,\ldots,s$ and
$\stackrel[l=1]{s}{\sum}m_{l}=r$ . We define $S$  a \emph{divisibility sequence} if  
 $$ S_{0}=0, S_{1}=1 \quad \forall n,m\geq 1, \quad S_{n}\mid S_{mn}. $$
 \end{definition}

Here we use the results on impulse sequences pointed out in previous section to retrieve the main result showed in \cite{BPP} with a different approach, based on determinants, giving a more detailed expression of the resultant sequence.
We recall the fundamental property which relates every recurrence sequence with suitable impulse sequences.
\begin{proposition}
	Every recurrence sequence $A=\left(a_n\right)_{n=0}^{+\infty}$ of order $r$ can be expressed in a unique way as a linear combination of $r$ impulse sequences of  order $r$ having the same minimal polynomial of $A$. More precisely we have
	\begin{equation}\label{uniquerepr}
	\forall n\geq 0 \quad a_{n}=\stackrel[k=0]{r-1}{\sum}a_{r-1-k}X_{n}^{(k)}
	\end{equation}
	where the terms $a_0,a_1,\ldots,a_{r-1}$, define the initial conditions of $A$ and, for $k=0,\ldots,r-1$, the recurrence sequences $(X_{n}^{(k)})_{n=0}^{+\infty}$ are the related impulse sequences.
	
\end{proposition}
\begin{proof}
	See, e.g. the fundamental book on recurrence sequences \cite{EPSW}.
\end{proof}

Now we are ready to prove the characterizing property of divisibility sequences pointed out in \cite{BPP}, in the general case of a divisibility sequence with minimal polynomial having multiple roots, giving a complete expression of the related resultant sequence.
\begin{theorem}\label{chardiv}
	Let $S$ be a non--degenerate recurrence sequence of order $r$ with minimal polynomial having $s$ distinct roots $\alpha_i$ with respective multiplicity $m_l$, $l=1,\ldots,s$ and  $\stackrel[l=1]{s}{\sum}m_{l}=r$. If $S$ is a divisibility sequence then  for all $n\geq 0$
	\begin{equation}\label{divrel}
	 S_{n}\mid D=n^{\stackrel[l=1]{s}{\sum}\left(\begin{array}{c}
		m_{l}\\
		2
		\end{array}\right)}\stackrel[l=1]{s}{\prod}\alpha_{l}^{\left(\begin{array}{c}
		m_{l}\\
		2
		\end{array}\right)(n-1)}\underset{1\leq i<j\leq s}{\prod}\left(\frac{\alpha_{j}^{n}-\alpha_{i}^{n}}{\alpha_{j}-\alpha_{i}}\right)^{m_{i}m_{j}}.
	\end{equation}
	
\end{theorem}
\begin{proof}
When $A=S$ from the equalities (\ref{uniquerepr}) we find that the following system of $r-1$ equations holds for every $n\geq1$
\begin{equation}
\begin{cases}
\stackrel[k=2]{r-1}{\sum}S_{r-1-k}X_{n}^{(k)}+S_{1}X_{n}^{(1)}+S_{0}X_{n}^{\left(0\right)}=S_n\\
\stackrel[k=2]{r-1}{\sum}S_{r-1-k}X_{2n}^{(k)}+S_{1}X_{2n}^{(1)}+S_{0}X_{2n}^{\left(0\right)}=S_{2n}\\
\cdots\\
\stackrel[k=2]{r-1}{\sum}S_{r-1-k}X_{(r-1)n}^{(k)}+S_{1}X_{(r-1)n}^{(1)}+S_{0}X_{(r-1)n}^{\left(0\right)}=S_{(r-1)n}
\end{cases}\label{eq:sistsn}
\end{equation}
If we consider $S$ as a divisibility sequence we have $S_0=0$ and $S_1=1$,
thus we can express $S_{1}$ using the
Cramer's rule applied to the coefficient matrix

\[
\left(X_{hn}^{(k)}\right)_{h=1,\ldots r-1,k=1,\ldots,r-1}
\]
whose determinant is $D\neq0$ as we proved in Theorem (\ref{detimp}).
We obtain
$$S_{1}=\frac{\det\left|\begin{array}{cccc}
		S_{n} & X_{n}^{(2)} & \cdots & X_{n}^{(r-1)}\\
		S_{2n} & X_{2n}^{(2)} & \cdots & X_{2n}^{(r-1)}\\
		\cdots & \cdots & \ddots & \cdots\\
		S_{(r-1)n} & X_{(r-1)n}^{(2)} & \cdots & X_{(r-1)n}^{(r-1)}
	\end{array}\right|}{D}=1$$
moreover $$S_n\mid \det\left|\begin{array}{cccc}
S_{n} & X_{n}^{(2)} & \cdots & X_{n}^{(r-1)}\\
S_{2n} & X_{2n}^{(2)} & \cdots & X_{2n}^{(r-1)}\\
\cdots & \cdots & \ddots & \cdots\\
S_{(r-1)n} & X_{(r-1)n}^{(2)} & \cdots & X_{(r-1)n}^{(r-1)}
\end{array}\right|$$ because $S_n$ divides all the entries $S_{hn}$, $h=1,\dots,r-1$, of the first column.
Therefore, observing that if $n=0$ $S_0|D=0$ , we clearly have
$$\forall n \geq0 \quad S_n \mid D.$$
\end{proof}

\begin{remark}

In particular, as a straightforward consequence of Theorem \ref{chardiv}, if the minimal polynomial of $S$ has all distinct roots, i.e.  we have $r=s$ and $m_{l}$=1 for all $l=1,\ldots,r$, equation (\ref{divrel}) becomes

\[
\forall n\geq0\quad S_{n}\mid\underset{1\leq i<j\leq r}{\prod}\left(\frac{\alpha_{j}^{n}-\alpha_{i}^{n}}{\alpha_{j}-\alpha_{i}}\right)
\]

since in this case obviously $\left(\begin{array}{c}
m_{l}\\
2
\end{array}\right)=0$  for every index $l$.
\end{remark}

\end{document}